\UseRawInputEncoding
\documentclass{amsart}

\usepackage{pstricks}
\usepackage{amssymb}
 \usepackage{pgfplots}

\usepackage{graphicx}
\usepackage{amsthm,amssymb,amsfonts,amsmath}
\usepackage[mathscr]{euscript}

\usepackage{bbm,color}
\usepackage{hyperref}
 \usepackage{mathptmx}      

\newcommand{\rme}{{\rm e}}
\newcommand{\rmd}{{\rm d}}
\newcommand{\la}{{\lambda}}
\newcommand{\ep}{{\varepsilon}}
\newcommand{\1}{\mathbbm{1}}

\newcommand{\sT}{{\mathcal T}}
\newcommand{\sP}{{\mathrm P}}
\newcommand{\sI}{{\mathcal I}}

\newtheorem{theo}{Theorem}[section]

\newtheorem{rem}{Remark}[section]

\newtheorem{cor}{Corollary}[section]

\begin{document}

\title{On telegraph processes, their first passage times and running extrema
}
\thanks{The work was supported by Russian Foundation for Basic Research (RFBR) 
and Chelyabinsk Region, project number 20-41-740020.}
\author{Nikita Ratanov\\
Chelyabinsk State University, 
129, Br.Kashirihykh, 
Chelyabinsk, 454001, Russia\\   
nikita.ratanov@csu.ru}
\date{}
\maketitle

\begin{abstract}
In this note, we present some ideas for describing the distributions of the 
running maximum/minimum, first passage times
and telegraphic meanders. Explicit formulae for joint distribution 
of the extrema, the number of velocity switches and the terminal position 
are derived using coupled integral equations technique.
\end{abstract}

\textbf{Keywords:} Asymmetric piecewise linear process; First passage time; Kac's scaling; 
Coupled integral equations; Telegraphic meander; Running maximum

\section{Some preliminaries}\label{sec1}

Let $\ep=\ep(t)\in\{0, 1\},\;t\ge0,$ be a two-state Markov process 
defined on the complete probability space $(\Omega, \mathcal F, \sP)$ and
governed by two alternating parameters $\la_0,\;\la_1,$
\begin{equation*}
\sP\{\ep(t+\rmd t)=i~|~\ep(t)=i\}=\exp(-\la_i\rmd t)+o(\rmd t),\qquad \rmd t\to0,\quad i\in\{0, 1\}.
\end{equation*}
Let  $\sP_i\{\cdot\}$ be the conditional probability $\sP\{\cdot~|~\ep(0)=i\}$
under the given initial state $\ep(0)=i.$

Denote by $\tau_1<\tau_2<\tau_3\ldots$ the switching times, $\tau_0=0.$
Let $N=N(t)$ be a right-continuous process that counts the switchings occurring up to time $t,$
 \begin{equation*}
N=N(t)=\inf\{n\ge0~|~\tau_n<t\},\qquad t>0,\qquad N(0)=0.
\end{equation*}

Consider the piecewise linear process 
$
\Gamma=\Gamma(t)=\int_0^t\gamma_{\ep(s)}\rmd s,\;t\ge0,
$ 
 which describes the evolution of the particle moving
with two velocities $\gamma_0,\;\gamma_1,\;\gamma_0>\gamma_1,$ 
alternating at random times $\tau_n,\;n\ge1.$ Process $\Gamma=\Gamma(t)$ 
can be regarded as an inhomogeneous asymmetric (integrated) telegraph process.
Throughout the paper, we use the following notations:
\begin{equation}
\label{def:xieta}
\xi_0=\xi_0(t, x)=\frac{x-\gamma_1t}{\gamma_0-\gamma_1},\qquad
\xi_1=t-\xi_0(t, x)=\frac{\gamma_0t-x}{\gamma_0-\gamma_1},
\end{equation}
\begin{equation}
\label{def:z}
z(t, x)=\la_0\la_1\xi_0(t, x)\xi_1(t, x),\qquad \gamma_1t<x<\gamma_0t,
\end{equation}
and
\begin{equation}
\label{def:thetapm}
\theta(t, x)
=\frac{\exp(-\la_0\xi_0(t, x)-\la_1\xi_1(t, x))}{\gamma_0-\gamma_1}.
\end{equation}
Notice that in the case of homogeneous telegraph process, that is
if $\la_0=\la_1=\la$ and $\gamma_0=-\gamma_1=\gamma,$ these notations are simplified to 
\[
\xi_0(t, x)=\frac{\gamma t+x}{2\gamma},\qquad \xi_1(t, x)=\frac{\gamma t-x}{2\gamma},\qquad
z(t, x)=\la^2(\gamma^2t^2-x^2),\qquad \theta(t, x)=\frac{\rme^{-\la t}}{2\gamma}.
\]
For instance, in these terms can be described the well-known joint distribution of $\Gamma(t)$ and 
$N(t),$ $P_i(\rmd x;\;t, x; n)=\sP_i\{\Gamma(t)\in\rmd x,\;N(t)=n\},$
$n\ge0,$ $i\in\{0, 1\}.$
When there is no velocity switching, the particle moves along the straight line, 
and the distribution of the particle's position is singular,
\begin{equation}
\label{eq:p00}
P_0(\rmd x; t, x; 0)=\rme^{-\la_0t}\delta_{\gamma_0t}(\rmd x),\qquad
P_1(\rmd x; t, x; 0)=\rme^{-\la_1t}\delta_{\gamma_1t}(\rmd x),
\end{equation}
where $\delta_a(\rmd x)$ is the (atomic) Dirac delta-measure. %
The probability density functions 
\[p_i(t, x;\;n)=P_i(\rmd x;\;t, x; n)/\rmd x,\qquad n\ge1,\]
can be written in terms of \eqref{def:xieta}-\eqref{def:thetapm}
separately for even and odd number of velocity switchings, 
see e.g. \cite[(4.1.10)-(4.1.11)]{KR}: 
\begin{equation}
\label{eq:p}
\begin{aligned}
    p_i(t, x;\;2n+1)&=\la_i\frac{z(t, x)^n}{n!^2} \cdot\theta(t, x)\1_{\{\gamma_1t<x<\gamma_0t\}}, \\
    p_i(t, x;\;2n+2)&=\la_0\la_1\xi_i(t, x)\frac{z(t, x)^n}{n!(n+1)!}\cdot\theta(t, x)\1_{\{\gamma_1t<x<\gamma_0t\}}, 
\end{aligned} i\in\{0, 1\}, \qquad n\ge0.
\end{equation}
Summing up, one can obtain the distribution of $\Gamma(t),\;t>0,$ 
\begin{equation}
\label{eq:P}
\begin{aligned}
&\sP_0\{\Gamma(t)\in\rmd x\}=\rme^{-\la_0t}\delta_{\gamma_0t}(\rmd x)
+\la_0\Big[\sI_0(t, x)\rmd x+\la_1\xi_0(t, x)\sI_1(t, x)\Big]\theta(t, x)\1_{\{\gamma_1t<x<\gamma_0t\}},\\
&\sP_1\{\Gamma(t)\in\rmd x\}=\rme^{-\la_1t}\delta_{\gamma_1t}(\rmd x)
+\la_1\Big[\sI_0(t, x)\rmd x+\la_0\xi_1(t, x)\sI_1(t, x)\Big]\theta(t, x)\1_{\{\gamma_1t<x<\gamma_0t\}}.
\end{aligned}
\end{equation}
Here functions $\sI_0$ and $\sI_1$ are defined by the series,
\begin{equation}
\label{def:sI}
 \sI_0(t, x)=\sum_{n=0}^\infty\frac{z(t, x)^n}{n!^2} =I_0\left(2\sqrt{z}\right)|_{z=z(t, x)},\qquad
\sI_1(t, x)=\sum_{n=0}^\infty\frac{z(t, x)^n}{n!(n+1)!} 
=\frac{I_1\left(2\sqrt{z}\right)}{\sqrt{z}}|_{z=z(t, x)},
\end{equation}
where $I_0$ and $I_1$ are the modified Bessel functions.

The main purpose of this note is to provide explicit formulae for the joint distribution of 
$\Gamma(t),$ the running maximum/minimum of $\Gamma,$ 
and the point at which the extremum is reached. 
We also introduce and study a telegraphic meander, 
which resembles the known brownian meanders, 
 see e.g. \cite{RevuzYor}.

Some preliminary results on the distribution of the maximum 
are known beginning with pioneering works by \cite{O90,pinsky,foong}. 
These results were later generalised to the inhomogeneous asymmetric case by 
 \cite{Luisa,BR,KR}.

Recently, a detailed analysis of a homogeneous case with drift was carried out in a series of papers by
\cite{cinque}. These results are heavily based on the fact that the arrival times $T_1,\ldots, T_n$
of the corresponding homogeneous Poisson process $N$ are uniformly distributed for a given number $N(t)=n.$
Unfortunately, the meticulous technique used in these papers 
does not work if $\la_0\ne\la_1.$ In this truly inhomogeneous case, a simple computation 
of conditional distributions under given number of velocity turns, $\{N(t)=n\},$ is unavailable,
 since arrivals of the counting process $N$ are distributed not uniformly. 
An explicit form of the conditional distribution of $\Gamma(t)$ for a given value of $N(t)$ can be found at
the end of Section 2 of \cite{JAP51}.

The text is organised as follows. First, we present a generalisation 
and detailing of known results on the distribution of first passage times
 (Section \ref{sec2}). 
Thereafter, the explicit representation of the meander distribution serves as an essential tool 
for obtaining the main result, Section \ref{sec3}.

\section{First passage of a telegraph process through a threshold}\label{sec2}
\setcounter{equation}{0}

Let $\sT(y)$  
be the time of the first passage by process $\Gamma=\Gamma(t)$ through the given threshold $y,$ 
\begin{equation*}
\sT(y)=\min\{t>0~|~\Gamma(t)=y\},
\end{equation*}
and $F_i(\rmd t, t, y)=\sP_i\{\sT(y)\in \rmd t\},\;i\in\{0, 1\},$ determine its distribution.

Let $F_i(\rmd t, t, y;\;n)=\sP_i\{\sT(y)\in\rmd t,\; N(t)=n\},\;n\ge0,\;i\in\{0, 1\},$ reflect the joint distribution
of the first passage time, $\sT(y),$  and the total number of velocity switchings, $N(t),$ 
occurring up to time $\sT(y)$.
Similarly to \eqref{eq:p00}, if there is no switching, then the distribution of $\sT(y)$ is singular, 
\begin{equation}
\label{eq:f00}
F_i\{\rmd t; t, y;\; 0\}=\rme^{-\la_it}\delta_{y/\gamma_i}(\rmd t),\quad i\in\{0, 1\}.
\end{equation}
As before, the probability density functions, $f_0(t, y;\; n)$ and $f_1(t, y;\; n),$ 
\[f_i(t, y;\; n)  =F_i(\rmd t; t, y;\;n)/\rmd t,   \qquad n\ge1,\qquad i\in\{0, 1\},\]
can be written separately
for different initial states 
and for an even and odd number of velocity changes. 
\begin{theo}\label{theo:fpt}
\begin{itemize}
  \item Let $y$ and both velocities  $\gamma_0,\;\gamma_1$ have the same sign\textup,
  and  $\Delta(y)$ be the support of the distribution of $\sT(y)$, that is, a segment with ends at the
 points $y/\gamma_0$ and $y/\gamma_1.$
   Then $\sT(y)$ is a. s. bounded\textup, $\sP\{\sT(y)\in\Delta(y)\}=1,$ and 
    \begin{equation}\label{eq:f0f1pp}
    \begin{aligned}
  f_i(t, y; 2n+1)  &= \la_i|\gamma_{1-i}|\frac{z(t, y)^n}{n!^2} \theta(t, y)\1_{\{t\in\Delta(y)\}},  \\
    f_i(t, y; 2n+2) &  =\la_0\la_1|\gamma_i|\xi_i(t, y)\frac{z(t, y)^{n}}{n!(n+1)!} \theta(t, y)\1_{\{t\in\Delta(y)\}},\\
 \end{aligned}i\in\{0, 1\},\qquad n\ge0. \end{equation}
 
 Further\textup,
     \begin{equation}\label{eq:Q01+}
\begin{aligned}
    F_0(\rmd t; t, y)=\rme^{-\la_0y/\gamma_0}&\delta_{y/\gamma_0}(\rmd t)\\
    +&\la_0\Big(|\gamma_1|\sI_0(t, y)+\la_1|\gamma_0|\xi_0(t, y)\sI_1(t, y)\Big)\theta(t, y)\1_{\{t\in\Delta(y)\}}\rmd t,   \\
 F_1(\rmd t; t, y)=\rme^{-\la_1y/\gamma_1}&\delta_{y/\gamma_1}(\rmd t)\\
    +&\la_1\Big(|\gamma_0|\sI_0(t, y)+\la_0|\gamma_1|\xi_1(t, y)\sI_1(t, y)\Big)\theta(t, y)\1_{\{t\in\Delta(y)\}}\rmd t, 
\end{aligned} 
\end{equation}

 In the cases $y<0<\gamma_1\le\gamma_0$ and $y>0>\gamma_0\ge\gamma_1,$ 
 the level $y$ is never reached\textup,
 $
 \sT(y)=\infty\;a.s.
$
  \item Let the velocities be of opposite signs\textup, $\gamma_0>0>\gamma_1,$ and $y>0.$
  Then\textup, we have 
  \[f_0(t, y;\;2n+1)\equiv0, \qquad f_1(t, y;\;2n+2)\equiv0\] 
  and 
   \begin{align}
     \label{eq:f12n+1}
     f_1(t, y;\;2n+1)=&\frac{\la_1}{\xi_0(t, y)}  \cdot\frac{z(t, y)^n\1_{\{t>y/\gamma_0\}}}{n!^2}\theta(t, y)
\cdot\left(y-\frac{\gamma_1}{n+1}\xi_1(t, y)\right),\\
   \label{eq:f02n}
     f_0(t, y;\;2n+2)=& 
     \la_0\la_1y\cdot\frac{z(t, y)^n}{n!(n+1)!}\theta(t, y)\1_{\{t>y/\gamma_0\}},\qquad n\ge0.
     \end{align}
\end{itemize}
Here $\xi_0(t, y),\;\xi_1(t, y),\;z(t, y)$ and $\theta(t, y)$
are defined by \eqref{def:xieta}-\eqref{def:thetapm}\textup.
  \end{theo}
 \begin{rem} 
 In the case of velocities with opposite signs\textup,  $\gamma_0>0>\gamma_1,$
 the formulae for $f_0(t, y;\; n)$ and $f_1(t, y;\; n),\;n\ge1,$ with negative threshold $y,$
 turn out to be symmetric to \eqref{eq:f12n+1}-\eqref{eq:f02n}\textup,
\begin{equation}
\label{eq:f01-}
f_0(t, y;\;n)|_{y<0}=f^{\leftrightarrow}_1(t, -y;\;n),\qquad
f_1(t, y;\;n)|_{y<0}=f^{\leftrightarrow}_0(t, -y;\;n),
\end{equation}
where $f^{\leftrightarrow}_0(t, \cdot;\;n)$ and $f^{\leftrightarrow}_1(t, \cdot;\;n)$ are determined by 
\eqref{eq:f12n+1}-\eqref{eq:f02n} with the following interchange of parameters\textup:
$
 \gamma_0\rightarrow-\gamma_1,\; \gamma_1\rightarrow-\gamma_0\;$ and $\la_0\leftrightarrow\la_1.
$
 Note that after these conversions\textup, we also have the interchange 
$\xi_0(t, y)\leftrightarrow\xi_1(t, -y),\;\xi_1(t, y)\leftrightarrow\xi_0(t, -y).$
Further\textup, summing up\textup, we obtain 
\begin{align}
\label{eq:Q0}
F_0(\rmd t; t, y)&=\begin{cases}
   \rme^{-\la_0y/\gamma_0}\delta_{y/\gamma_0}(\rmd t)
   +\la_0\la_1y\sI_1(t, y)\theta(t, y)\1_{\{t>y/\gamma_0\}}\rmd t,   & y>0, \\ \\
\dfrac{\la_0}{\xi_1(t, y)} \left(-y\sI_0(t, y)+\gamma_0\xi_0(t, y)\sI_1(t, y)\right)\theta(t, y)\1_{\{t>y/\gamma_1\}}\rmd t,
 & y<0.
\end{cases}
 \\
 \label{eq:Q1}
    F_1(\rmd t; t, y)&=\begin{cases}
     \dfrac{\la_1}{\xi_0(t, y)}\Big(y\sI_0(t, y)-\gamma_1\xi_1(t, y)\sI_1(t, y)\Big)
     \theta(t, y)\1_{\{t>y/\gamma_0\}}\rmd t,&
 y>0, \\ \\
    \rme^{-\la_1y/\gamma_1}\delta_{y/\gamma_1}(\rmd t)
 -\la_0\la_1y\sI_1(t, y)\theta(t, y)\1_{\{t>y/\gamma_1\}},     & y<0,
\end{cases}  
\end{align}
cf. \textup{\cite[Theorem 3.1]{JAP51}.}
  \end{rem}
\begin{proof}
First, let both velocities be positive, $\gamma_0>\gamma_1>0$. 
In this case, 
$\sT(y),\;y>0,$ is a. s. bounded, process $\Gamma(t),\;t\ge0,$ is a subordinator, and, hence, 
$f_i(t, y;\;n)\rmd t\sim p_i(t, y;\; n)\rmd y,\; \rmd t\to0,
$
where $\rmd y$ is the increment of $\Gamma(t)$ corresponding to the time increment $\rmd t.$
Therefore, $f_i(t, y;\;n)=\gamma_{\ep(t)}p_i(t, y;\; n),$ and
formulae \eqref{eq:f0f1pp}, $\gamma_0>\gamma_1>0,\;y>0,$ follow from \eqref{eq:p}. 
The case of both negative velocities is symmetric. 
Expressions \eqref{eq:Q01+}  for $F_0(\rmd t; t, y)$ and $F_1(\rmd t; t, y)$
follow by summing the formulae \eqref{eq:f0f1pp}.

Let the velocities have opposite signs, $\gamma_0>0>\gamma_1,$ and the threshold is positive, $y>0.$	
Since the first crossing of a positive threshold always occurs at the positive velocity, we have
\[F_1(\rmd t; t, y;\;0)|_{y>0}\equiv0,\qquad f_1(t, y;\;2n+2)|_{y>0}\equiv0,\qquad f_0(t, y;\;2n+1)|_{y>0}\equiv0, \quad n\ge0,\quad y>0.\]
Further, note that starting from the state $0=\ep(0),$ the particle first crosses the positive level $y$ only
after an even number of switchings, and the first turn occurs before the
time $y/\gamma_0,$  that is before the particle reaches the threshold $y.$
Similarly, if the particle starts from $1=\ep(0),$ it must perform an odd number of switchings,  
and the first turn must be before the time $\xi_1(t, y).$
Conditioning on the first velocity switching, 
we obtain the following sequence of coupled integral equations, for $n\ge1,$
\begin{equation}
\label{eq:f0f1int}
\left\{
\begin{aligned}
  f_0(t, y;\;2n) &
  =\int_0^{y/\gamma_0}\la_0\rme^{-\la_0\tau}f_1(t-\tau, y-\gamma_0\tau;\;2n-1)\rmd\tau,  
 \\
 f_1(t, y;\;2n+1) &
    =\int_0^{\xi_1(t, y)}\la_1\rme^{-\la_1\tau}f_0(t-\tau, y-\gamma_1\tau;\;2n)\rmd\tau,
\end{aligned}\qquad y>0.
\right.
\end{equation}

Equations \eqref{eq:f0f1int} can be solved explicitly. 
For example, for  $n=0,$ the second equation of this system, due to \eqref{eq:f00}, takes the form
\begin{equation}
\label{eq:f11+}\begin{aligned}
f_1(t, y;\;1)=&\int_0^{t}\la_1\rme^{-\la_1\tau}F_0(\rmd \tau;\;t-\tau, y-\gamma_1\tau;\;0)\\
=&\frac{\la_1\gamma_0}{\gamma_0-\gamma_1}\rme^{-\la_0\xi_0(t, y)-\la_1\xi_1(t, y)}
=\la_1\gamma_0\theta(t, y)\1_{\{t>y/\gamma_0\}}.
\end{aligned}
\end{equation}

We continue working with system \eqref{eq:f0f1int} using the identities 
\begin{align}
\label{eq:xieta0}
   \xi_0(t-\tau, y-\gamma_0\tau)\equiv\xi_0(t, y)-\tau, &\qquad\xi_1(t-\tau, y-\gamma_0\tau) \equiv\xi_1(t, y),
  \\
    \label{eq:xieta1}
 \xi_0(t-\tau, y-\gamma_1\tau)\equiv\xi_0(t, y), &   \qquad\xi_1(t-\tau, y-\gamma_1\tau) \equiv\xi_1(t, y)-\tau
\end{align}
and
\begin{equation}
\label{eq:y}
\gamma_0\xi_0(t, y)+\gamma_1\xi_1(t, y)\equiv y,
\end{equation}
which are obvious by definition, \eqref{def:xieta}. 

For $n=0,$ formula \eqref{eq:f12n+1} is proved, see \eqref{eq:f11+} and \eqref{eq:y}.
The subsequent formulae in \eqref{eq:f02n} and \eqref{eq:f12n+1} follow by induction.

Substituting  \eqref{eq:f12n+1} (with $n-1$ instead of $n$)
  into the first equation of \eqref{eq:f0f1int}, due to \eqref{eq:xieta0} and \eqref{eq:y}
we obtain
 \[
\begin{aligned}
   f_0(t, y; 2n) & 
   =\int_0^{y/\gamma_0} \la_0\frac{\la_1}{\xi_0-\tau}
\left(y-\gamma_0\tau-\frac{\gamma_1}{n}\xi_1\right)  
 \frac{(\la_0\la_1)^{n-1}(\xi_0-\tau)^{n-1}\xi_1^{n-1}}{(n-1)!^2}\rmd\tau\cdot\theta(t, y)\1_{\{t>y/\gamma_0\}}\\
  & =\frac{(\la_0\la_1)^{n}\xi_1^{n-1}}{(n-1)!^2}\int_0^{y/\gamma_0}
\left(\gamma_0(\xi_0-\tau)+\frac{n-1}{n}\gamma_1\xi_1\right)
(\xi_0-\tau)^{n-2}\rmd\tau\cdot\theta(t, y)\1_{\{t>y/\gamma_0\}}\\
&=
 \frac{(\la_0\la_1)^{n}\xi_1^{n-1}}{(n-1)!^2} 
\int_{-\gamma_1\xi_1/\gamma_0}^{\xi_0}\left(\gamma_0u^{n-1}
+\frac{n-1}{n}\gamma_1\xi_1 u^{n-2}\right)\rmd u
\cdot \theta(t, y)\1_{\{t>y/\gamma_0\}},
\end{aligned} 
 \]
where $\xi_0=\xi_0(t, y)$ and $\xi_1=\xi_1(t, y)$.  After integration, by virtue of \eqref{eq:y}, we obtain
 \[
    f_0(t, y; 2n) =
\frac{(\la_0\la_1)^{n}\xi_1^{n-1}}{(n-1)!^2} \cdot
   \frac{\gamma_0\xi_0^{n}+\gamma_1\xi_1\xi_0^{n-1}}{n}\cdot \theta(t, y)
    =y\frac{(\la_0\la_1)^{n}(\xi_0\xi_1)^{n-1}}{(n-1)!n!}\cdot\theta(t, y),
    \qquad t>y/\gamma_0,
\]
which confirms \eqref{eq:f02n}.

Similarly, substituting  \eqref{eq:f02n} (with $n-1$ instead of $n$) 
into the second equation of \eqref{eq:f0f1int}, due to \eqref{eq:xieta1} and \eqref{eq:y},
 we get
 \[
\begin{aligned}
f_1(t, n;\;2n+1)&=\la_0\la_1\int_0^{\xi_1}\la_1(y-\gamma_1\tau)
\frac{\left[\la_0\la_1\xi_0(\xi_1-\tau)\right]^{n-1}}{(n-1)!n!}\rmd\tau\cdot\theta(t, y)\1_{\{t>y/\gamma_0\}}\\
&=\frac{\la_0^n\la_1^{n+1}\xi_0^{n-1}}{(n-1)!n!}\int_0^{\xi_1}
\left(\gamma_0\xi_0+\gamma_1(\xi_1-\tau)\right)(\xi_1-\tau)^{n-1}\rmd\tau\cdot\theta(t, y)\1_{\{t>y/\gamma_0\}}
\\
&=\frac{\la_0^n\la_1^{n+1}\xi_0^{n-1}}{(n-1)!n!}
\left(\gamma_0\xi_0\frac{\xi_1^n}{n}+\gamma_1\frac{\xi_1^{n+1}}{n+1}\right)\cdot\theta(t, y)\1_{\{t>y/\gamma_0\}}\\
&=\frac{\la_0^n\la_1^{n+1}\xi_0^{n-1}\xi_1^n}{n!^2}
\left(y-\frac{\gamma_1}{n+1}\xi_1\right)\cdot\theta(t, y)\1_{\{t>y/\gamma_0\}},
\end{aligned} 
\]
which coincides with \eqref{eq:f12n+1}.
\end{proof}
  
  \begin{rem}
  The idea of using coupled integral equations of the form \eqref{eq:f0f1int} 
to analyse the properties of telegraph-like processes
(instead of the classic technique which is based on the 
differential equations) occasionally appears in literature, see e.g.\textup{ \cite[Lemma 5.1]{DiC2001},
\cite[formula (5.6)]{SZ2004}, \cite[Chapter 5]{zacks}, \cite[(4.1.2)]{KR}.}

Formulae \eqref{eq:Q0}-\eqref{eq:Q1} were previously derived by a slightly different method, see 
\textup{\cite{BR,JAP51}.}
At first glance, it seems that formulae \eqref{eq:f02n}-\eqref{eq:f12n+1} can be obtained  by expanding
the Bessel functions in formulae \eqref{eq:Q0}-\eqref{eq:Q1}. However, to prove this, we need to show that 
the $n$-th term of these series must be equal to $f_\cdot(t, y;\; n),$ which is not so evident.
\end{rem}

\begin{rem}\label{rem:hatphi}
The distribution of the time to first reach of the threshold $y,\;y>0,$ with a simultaneous velocity reversal
at this time follows from
\begin{equation}
\label{def:hatfi}
\sP_i\left\{\sT(y)\in\rmd t,\;N(t)=n, N(t+)=n+1\right\}=
     \la_0\cdot\sP_i\left\{\sT(y)\in\rmd t,\;N(t)=n\right\}, 
     \end{equation}
   $y>0,\;  i\in\{0, 1\},\; n\ge0.$
For $n=0,$ these equalities hold\textup, because $\sP_1\left\{\sT(y)\in \rmd t,\;N(t)=0\right\}\equiv0,$ $t>0,$ and
\[\sP_0\left\{\sT(y)\in\rmd t,\;N(t)=0, N(t+)=1\right\}
=\sP_0\{\tau_1\in\rmd t\text{ and }y=\gamma_0t\}=
\la_0\rme^{-\la_0t}\delta_{y/\gamma_0}(\rmd t),\] which is $\la_0F_0(\rmd t, t, y;\;0).$ 

For $i=0$ and even $n,\; i=1$ and odd $n,$
 equality \eqref{def:hatfi} follows in view of the integral equations \eqref{eq:f0f1int}.
 Equalities \eqref{def:hatfi} with $i=0$ and odd $n$ $(i=1$ and even $n)$ are trivial\textup, $0=0.$
 \end{rem}

It is interesting to look at formulae \eqref{eq:Q0}-\eqref{eq:Q1}
in light of the Kac's rescaling. Exactly, let $\la_0,\la_1\to\infty$ and
\begin{equation}
\label{eq:la0la1}
\frac{\la_0}{\la_1}\to\nu^2,\qquad \nu>0.
\end{equation}
The Kac condition is assumed for the two states separately: let  
 $\gamma_0\to+\infty,\;\gamma_1\to-\infty$ and 
\begin{equation}
\label{eq:gamma01}
\frac{\gamma_0}{\sqrt{\la_0}}\to\sigma_0,\qquad
\frac{\gamma_1}{\sqrt{\la_1}}\to-\sigma_1,
\end{equation}
where $\sigma_0,\;\sigma_1>0.$ Finally, we assume that
\begin{equation}
\label{eq:delta}
\frac{\gamma_0\la_1+\gamma_1\la_0}{\la_0+\la_1}\to\delta.
\end{equation}

It is known that under conditions \eqref{eq:la0la1}-\eqref{eq:delta}, 
process $\Gamma(t)$ weakly converges to the scaled Brownian motion with drift, $\Sigma \cdot W(t)+\delta t,\;t>0,$
where 
\begin{equation}
\label{def:Sigma}
\Sigma=\dfrac{\sigma_0\sigma_1}{\sqrt{(\sigma_0^2+\sigma_1^2)/2}}. 
\end{equation}
See \cite[Section 5]{JAP51} for detailed definitions and comments.
In the symmetric case, $\la_0=\la_1$ and $\gamma_0=-\gamma_1,$ we have $\delta=0$ and $\sigma_0=\sigma_1=\Sigma$,
so that the telegraph process $\Gamma(t)$ weakly converges to $\Sigma \cdot W(t),\;t>0.$

The following result seems natural: 
under the Kac scaling \eqref{eq:la0la1}-\eqref{eq:delta}, the distribution of $\sT(y)$
converges to the first passage time distribution of Brownian motion.  
\begin{cor}[Cf  \cite{KS}]\label{cor1}
Under the Kac scaling defined by \eqref{eq:la0la1}-\eqref{eq:delta}\textup,
\begin{equation}
\label{eq:limF0F1}
F_0(\rmd t; t, y),\;F_1(\rmd t; t, y)\to
\frac{y}{\sqrt{2\pi}\Sigma t^{3/2}}\exp\left(-\frac{(y-\delta t)^2}{2\Sigma^2t}\right),
\end{equation}
where $\Sigma$ is given by \eqref{def:Sigma}.
\end{cor}
\begin{proof}
We apply the scaling conditions \eqref{eq:la0la1}-\eqref{eq:delta} to \eqref{eq:Q0}-\eqref{eq:Q1} and use 
the asymptotic expansion of the Bessel functions
$
I_k(x)\sim\dfrac{\exp(x)}{\sqrt{2\pi x}},\; x\to\infty,
$
see e.g. \cite[8.451]{GR}.
\end{proof}

\section{Telegraphic meanders and running extrema}\label{sec3}
\setcounter{equation}{0}
We define the sets of always negative and always positive trajectories with a given number of switchings, 
say a \emph{negative telegraphic meander},
\begin{equation}
\label{def:meander-}
\mathfrak m^-(t, x;\;n)=\{\Gamma(t)\in\rmd x,\;M_t=0,\; N(t)=n\},
\end{equation}
and a \emph{positive telegraphic meander},
\begin{equation}
\label{def:meander+}
\mathfrak m^+(t, x;\;n)=\{\Gamma(t)\in\rmd x,\;m_t=0,\; N(t)=n\}.
\end{equation}
where $m_t:=\min_{u\in[0, t]}\Gamma(u)$ and $M_t:=\max_{u\in[0, t]}\Gamma(u)$ 
are running minimum and running maximum. 

If both velocities are of the same sign, then $\Gamma(t),\;t>0,$  preserves the sign on an arbitrary time interval.
In the case of both negative velocities, $M_t=0$ and $m_t<0, \;\forall t>0.$
In this case, 
\[
\sP_i\{\mathfrak m^-(t, x;\;n)\}=
\sP_i\{\Gamma(t)\in\rmd x,\;N(t)=n\}
\qquad\text{and}\qquad \sP_i\{\mathfrak m^+(t, x;\; n)\}\equiv0,
\quad i\in\{0, 1\};
\] 
 if the velocities are positive, we have 
 \[
 \sP_i\{\mathfrak m^-(t, x;\; n)\}\equiv0,\qquad\text{and}\qquad\sP_i\{\mathfrak m^+(t, x;\;n)\}=
\sP_i\{\Gamma(t)\in\rmd x,\;N(t)=n\},
\quad i\in\{0, 1\},
 \]
 the explicit formulae for $\sP_i\{\Gamma(t)\in\rmd x,\;N(t)=n\}$ are given by
 \eqref{eq:p00}-\eqref{eq:p}.

 Let the velocities be of opposite signs, $\gamma_0>0>\gamma_1.$ 
Therefore,  $M_t|_{\ep(0)=0}>0,\;m_t|_{\ep(0)=1}<0$ and
$
\sP_0\{\mathfrak m^-(t, x;\;n)\}\equiv0,\;\sP_1\{\mathfrak m^+(t, x;\;n)\}\equiv0.
$
We are interested in explicit expressions for the distributions of negative and positive meanders,
\begin{equation}
\label{def:Gn}
G^-(\rmd x, t, x;\; n)=\sP_1\{\mathfrak m^-(t, x;\;n)\},\quad x<0;\quad
G^+(\rmd x, t, x;\; n)=\sP_0\{\mathfrak m^+(t, x;\;n)\},\quad x>0.
\end{equation}
For brevity, we focus on the positive meander, i. e. on $G^+.$ 
Note that in the absence of switches, i.e. $N(t)=0,$ 
 the distribution of $\Gamma(t)$ is singular, see \eqref{eq:p00},
\begin{equation}
\label{eq:g-0g+0}
G^+(\rmd x, t, x; 0)
    =\rme^{-\la_0t}\delta_{\gamma_0t}(\rmd x),\qquad t>0.
\end{equation}
Consider the probability density functions, $g^+(t, x;\; n)=\sP_0\{\mathfrak m^+(t, x; n)\}/\rmd x,\; n\ge1.$


\begin{theo}\label{theo:alwayspositive}
Let $\gamma_0>0>\gamma_1.$ Functions 
$g^+(t, x;\; n),\;x>0,$ are specified explicitly\textup, separately for odd and even $n,$ by the expressions
    \begin{align}
     \label{eq:g+2n+1}
    g^+(t, x; 2n+1)=&\frac{\la_0}{\xi_0(t, x)} \frac{z(t, x)^n\1_{\{0<x<\gamma_0t\}}}{n!^2}
     \theta(t, x) \cdot\left(x-\frac{\gamma_1}{n+1}\xi_1(t, x)\right)/\gamma_0,\\  
  \label{eq:g+2n}
    g^+(t, x; 2n+2)=&\la_0\la_1\frac{z(t, x)^{n}\1_{\{0<x<\gamma_0t\}}}{n!(n+1)!}
    \theta(t, x)\cdot x/\gamma_0,\qquad n\ge0. 
\end{align}
\end{theo}
\begin{rem}
By summing up \eqref{eq:g-0g+0}-\eqref{eq:g+2n}\textup, we obtain
\begin{equation}
\label{eq:QQ+}\begin{aligned}
G^+(\rmd x, t, x)&=\sP_0\{\Gamma(t)\in\rmd x,\;m_t=0\}\\&=\rme^{-\la_0t}\delta_{\gamma_0t}(\rmd x)
+\frac{\la_0}{\gamma_0}\left[\frac{x}{\xi_0(t, x)}\sI_0(t, x)
+\left(\la_1x-\gamma_1\frac{\xi_1(t, x)}{\xi_0(t, x)}\right)\sI_1(t, x)\right]\theta(t, x)\rmd x.
\end{aligned}\end{equation}

The distribution of the negative meander follows by symmetry\textup:
the formulae for $g^-(t, x;\; n),$ $\gamma_1t<x<0,$ have the form \eqref{eq:g+2n+1}-\eqref{eq:g+2n}
with interchange $0\leftrightarrow1,$ and after summing\textup,
\begin{equation}
\label{eq:QQ-}\begin{aligned}
G^-(\rmd x, t, x)&=\sP_1\{\Gamma(t)\in\rmd x,\;M_t=0\}\\&=\rme^{-\la_1t}\delta_{\gamma_1t}(\rmd x)
+\frac{\la_1}{\gamma_1}\left[\frac{x}{\xi_1(t, x)}\sI_0(t, x)
+\left(\la_0x-\gamma_0\frac{\xi_0(t, x)}{\xi_1(t, x)}\right)\sI_1(t, x)\right]\theta(t, x)\rmd x.
\end{aligned}\end{equation}
\end{rem}

\begin{proof} 
The proof of the theorem is based on the following observation:
each meander path $\gamma(s, x),$ $s\in[0, t],$ 
considered in reverse time, that is, $\gamma(t-s, x),\;s\in[0, t],$
 matches the trajectory starting at $x$ with first passage through the origin at time $t$.
  Therefore, the probability density functions
$g^+(t, x;\; n)$ can be obtained similarly to \eqref{eq:f12n+1}-\eqref{eq:f02n}. 

In contrast to the proof of Theorem \ref{theo:fpt}, the coupled integral equations 
for $g^+(t, x;\;\cdot)$ are written out by
conditioning on the \emph{last} velocity change, cf Remark \ref{rem:hatphi}. We have
\begin{equation}
\label{eq:g+eq}
\left\{
\begin{aligned}
   g^+(t, x;2n+2)&  =\int_0^{x/\gamma_0}\la_1\rme^{-\la_0s}g^+(t-s, x-\gamma_0s;2n+1)\rmd s,\\
    g^+(t, x, 2n+1)& = \int_0^{\xi_1(t, x)}\la_0\rme^{-\la_1s}g^+(t-s, x-\gamma_1s; 2n)\rmd s.  
\end{aligned}
\right.
\end{equation}
%
Formulae \eqref{eq:g+2n+1}-\eqref{eq:g+2n}  for always positive paths can be proved by induction 
similarly to the proof of Theorem \ref{theo:fpt}. 
Formula \eqref{eq:QQ+} follows by summing up \eqref{eq:g+2n+1}-\eqref{eq:g+2n}.
\end{proof}

The above preparation allows to receive the joint distribution of the running extrema $m_t\;(M_t),$
 the time $\zeta_t^m\; (\zeta_t^M)$ to reach the running extremum, and the terminal position $\Gamma(t),$
 reached after $N(t)=n$ velocity switchings.
 
If both velocities are positive, $\gamma_0>\gamma_1>0,$ then $\zeta_t^m=0,\; m_t=0$ 
and $\zeta_t^M=t,\;M_t=\Gamma(t)$ a. s. Therefore, in this case, for $i\in\{0, 1\},$
 \begin{equation*}
 \begin{aligned}
\sP_i\left\{\zeta_t^m\in\rmd s,\;m_t\in\rmd y,\;\Gamma(t)\in\rmd x,\; N(t)=n\right\}
=&\delta_0(\rmd s)\delta_0(\rmd y)p_i(t, x;\;n)\rmd x,\\%
\sP_i\left\{\zeta_t^M\in\rmd s,\;M_t\in\rmd y,\;\Gamma(t)\in\rmd x,\; N(t)=n\right\}
=&\delta_t(\rmd s)\delta_{x}(\rmd y)p_i(t, x;\;n)\rmd x,
\end{aligned}\qquad \gamma_1t<x<\gamma_0t;
\end{equation*}
similarly, for $0>\gamma_0>\gamma_1,$
 \begin{equation*}
 \begin{aligned}
\sP_i\left\{\zeta_t^m\in\rmd s,\;m_t\in\rmd y,\;\Gamma(t)\in\rmd x,\; N(t)=n\right\}
=&\delta_t(\rmd s)\delta_{x}(\rmd y)p_i(t, x;\;n)\rmd x,\\
\sP_i\left\{\zeta_t^M\in\rmd s,\;M_t\in\rmd y,\; \Gamma(t)\in\rmd x,\; N(t)=n\right\}
=&\delta_0(\rmd s)\delta_{0}(\rmd y)p_i(t, x;\;n)\rmd x,
\end{aligned}\qquad \gamma_1t<x<\gamma_0t,
\end{equation*} 
the formulae for $p_i(t, x;\;n)$ are given by \eqref{eq:p}. 
In these two cases, 
the distributions of $\langle\zeta_t^m,\;m_t,\;\Gamma(t)\rangle$ and  $\langle\zeta_t^M,\;M_t,\;\Gamma(t)\rangle$
can be obtained by summing up, see  \eqref{eq:P}.
 
Let the velocities be of opposite signs, $\gamma_0>0>\gamma_1.$
The singular components of these distributions
corresponding to $\{\zeta_t^m=0\}$ and $\{\zeta_t^M=0\}$ are represented as follows.
Each time the running minimum or running maximum is zero, 
$m_t=0$ or $M_t=0$, the trajectory follows the meander, and
 \begin{equation}
\label{eq:zeta0m}
 \begin{aligned}
\sP_0\left\{\zeta_t^m=0,\;m_t=0,\;\Gamma(t)\in\rmd x,\;N(t)=n\right\}=&G^+(\rmd x; t, x;\;n),\qquad x>0,\\
\sP_1\left\{\zeta_t^M=0,\;M_t=0,\;\Gamma(t)\in\rmd x,\;N(t)=n\right\}=&G^-(\rmd x; t, x;\; n),\qquad x<0,
\end{aligned}
\end{equation}
see \eqref{def:Gn}-\eqref{eq:g+2n}. Therefore,
\begin{equation}
\label{eq:P0P1meander}
\begin{aligned}
 \sP_0\left\{\zeta_t^m=0,\;m_t=0,\;\Gamma(t)\in\rmd x\right\}=&G^+(\rmd x;\;t, x),   \\
\sP_1\left\{\zeta_t^M=0,\;M_t=0,\;\Gamma(t)\in\rmd x\right\}=&G^-(\rmd x;\; t, x),
\end{aligned}
\end{equation}
see \eqref{eq:QQ+}-\eqref{eq:QQ-}.

Further, since $\Gamma(t)$ arrives at the minimum, $m_t$ (maximum, $M_t$)  
with the negative (positive) velocity, 
the negative running minimum $x$ is reached at the last time, that is $\zeta_t^m=t,$ with probabilities   
 \begin{equation}
\label{eq:zetat0m}
 \begin{aligned}
\sP_0\left\{\zeta_t^m=t,\;m_t=\Gamma(t)\in\rmd x,\;N(t)=n\right\}=\begin{cases}
      0,& \text{if $n$ is even }, \\
      f_0(t, x;\;n)\rmd x,& \text{if $n$ is odd},
\end{cases}\\
\sP_1\left\{\zeta_t^m=t,\;m_t=\Gamma(t)\in\rmd x,\;N(t)=n\right\}=\begin{cases}
      0,& \text{if $n$ is odd }, \\
      f_1(t, x;\;n)\rmd x,& \text{if $n$ is even},
\end{cases}
\end{aligned}\qquad x<0;
\end{equation}
and, respectively, for $\zeta^M_t=t,$
\begin{equation}
\label{eq:zetat0M}
 \begin{aligned}
\sP_0\left\{\zeta_t^M=t,\;M_t=\Gamma(t)\in\rmd x,\;N(t)=n\right\}=\begin{cases}
      0,& \text{if $n$ is odd }, \\
      f_0(t, x;\;n)\rmd x& \text{if $n$ is even},
\end{cases}\\
\sP_1\left\{\zeta_t^M=t,\;M_t=\Gamma(t)\in\rmd x,\;N(t)=n\right\}=\begin{cases}
      0,& \text{if $n$ is even }, \\
      f_1(t, x;\;n)\rmd x& \text{if $n$ is odd},
\end{cases}
\end{aligned}\qquad x>0,
\end{equation}
where $f_0(t, x;\; n)$ and $f_1(t, x;\; n)$ are given by \eqref{eq:f00}, \eqref{eq:f12n+1}-\eqref{eq:f01-}.
Summing up, we obtain
\begin{equation}
\label{eq:Pm}\begin{aligned}
\sP_0\left\{\zeta_t^m=t,\;m_t=\Gamma(t)\in\rmd x\right\}=&
\frac{\la_0}{\xi_1(t, x)}\left(-x\sI_0(t, x)+\gamma_0\xi_0(t, x)\sI_1(t, x)\right)
\theta(t, x)\1_{\{\gamma_1t<x<0\}}\frac{\rmd x}{-\gamma_1},\\
\sP_1\left\{\zeta_t^m=t,\;m_t=\Gamma(t)\in\rmd x\right\}=&
\left(\rme^{-\la_1x/\gamma_1}\delta_{\gamma_1t}(\rmd x)-\la_0\la_1x\sI_1(t, x)\right)
\theta(t, x)\1_{\{\gamma_1t<x<0\}}\frac{\rmd x}{-\gamma_1},\\
\sP_0\left\{\zeta_t^M=t,\;M_t=\Gamma(t)\in\rmd x\right\}=&
\left(\rme^{-\la_0x/\gamma_0}\delta_{\gamma_0t}(\rmd x)+\la_0\la_1x\sI_1(t, x)\right)
\theta(t, x)\1_{\{0<x<\gamma_0t\}}\frac{\rmd x}{\gamma_0},\\
\sP_1\left\{\zeta_t^M=t,\;M_t=\Gamma(t)\in\rmd x\right\}=&
\frac{\la_1}{\xi_0(t, x)}\left(x\sI_0(t, x)-\gamma_1\xi_1(t, x)\sI_1(t, x)\right)
\theta(t, x)\1_{\{0<x<\gamma_0t\}}\frac{\rmd x}{\gamma_0}.
\end{aligned}\end{equation}
Compare formulae \eqref{eq:zetat0m}-\eqref{eq:Pm}
with similar formulae obtained in the symmetric case, see \cite{cinque}.

The ``regular" component of the distribution is obtained 
by weighing the compound paths consisting of a passage to a minimum (maximum) value, 
switching the velocity at this moment, 
then moving along the meander above (below) the reached level in the remaining time.
The corresponding probabilities for $0<s<t$ turn out to be
\begin{equation*}
\begin{aligned}
\sP_0\left\{\zeta_t^m\in\rmd s,\;m_t\in\rmd y,\;\Gamma(t)\in\rmd x,\;N(t)=n\right\}&\\
=
\la_1\sum_{1\le k\le [n/2]} F_0(\rmd s, s, y;\;2k-1)&G^+(\rmd x, t-s, x-y;\;n-2k)\rmd y,\qquad n\ge2,
\\
\sP_1\left\{\zeta_t^m\in\rmd s,\;m_t\in\rmd y,\;\Gamma(t)\in\rmd x,\;N(t)=n\right\}&\\
=
\la_1\sum_{0\le k\le [n/2]} F_1(\rmd s, s, y;\;2k)&G^+(\rmd x, t-s, x-y;\;n-2k-1)\rmd y,\qquad n\ge1,
\end{aligned}\end{equation*}
\[0<s<t,\qquad y<0\wedge x,\]
and 
\begin{equation*}
\begin{aligned}
\sP_0\left\{\zeta_t^M\in\rmd s,\;M_t\in\rmd y,\;\Gamma(t)\in\rmd x,\;N(t)=n\right\}&\\
=
\la_0\sum_{0\le k\le [(n-1)/2]} F_0(\rmd s, s, y;\;2k)&G^-(\rmd x, t-s, x-y;\;n-2k-1)\rmd y,\qquad n\ge1,
\\
\sP_1\left\{\zeta_t^M\in\rmd s,\;M_t\in\rmd y,\;\Gamma(t)\in\rmd x,\;N(t)=n\right\}&\\
=
\la_0\sum_{1\le k\le [n/2]} F_1(\rmd s, s, y;\;2k-1)&G^-(\rmd x, t-s, x-y;\;n-2k)\rmd y,\qquad n\ge2,
\end{aligned}\end{equation*}
\[0<s<t,\qquad y>0\vee x.\]
Here $F_0(\rmd s;\; \cdot, \cdot;\; \cdot)$ and $F_1(\rmd s;\; \cdot, \cdot;\; \cdot)$
are written out in \eqref{eq:f00},  \eqref{eq:f12n+1}-\eqref{eq:f01-}, 
and $G^\mp(\rmd x;\;\cdot, \cdot;\;\cdot)$ are determined by formulae
\eqref{def:Gn}-\eqref{eq:g+2n}.

By summing up, we obtain the distributions of  
$\langle\zeta_t^m,\;m_t,\;\Gamma(t)\rangle$ and  $\langle\zeta_t^M,\;M_t,\;\Gamma(t)\rangle$
in terms of sum of Bessel functions:
\[
\begin{aligned}
\sP_0\left\{\zeta_t^m\in\rmd s,\;m_t\in\rmd y,\;\Gamma(t)\in\rmd x\right\}& 
=\delta_0(\rmd s)\delta_0(\rmd y)G^+(\rmd x, t, x)\1_{\{x>0\}}  \\
+&\frac{\la_0\delta_t(\rmd s)\rmd x}{-\gamma_1\xi_1(t, x)}\left(-x\sI_0(t, x)
+\gamma_0\xi_0(t, x)\sI_1(t, x)\right)\1_{\{t>x/\gamma_1,\;x<0\}}\delta_x(\rmd y)\\
    +&  \la_1F_0(\rmd s, s, y)G^+(\rmd x, t-s, x-y)\rmd y,\\
    \sP_1\left\{\zeta_t^m\in\rmd s,\;m_t\in\rmd y,\;\Gamma(t)\in\rmd x\right\}& 
=\left(\rme^{-\la_1x/\gamma_1}\delta_{\gamma_1t}(\rmd x)-\la_0\la_1x\sI_1(t, x)\right)
\theta(t, x)\1_{\{\gamma_1t<x<0\}}\frac{\rmd x}{-\gamma_1}\\
&+\la_1F_1(\rmd s, s, y)G^+(\rmd x, t-s, x-y)\rmd y,
    \\
\sP_0\left\{\zeta_t^M\in\rmd s,\;M_t\in\rmd y,\;\Gamma(t)\in\rmd x\right\}& 
=\left(\rme^{-\la_0x/\gamma_0}\delta_{\gamma_0t}(\rmd x)+\la_0\la_1x\sI_1(t, x)\right)
\theta(t, x)\1_{\{0<x<\gamma_0t\}}\frac{\rmd x}{\gamma_0}\\
&+\la_0F_0(\rmd s, s, y)G^-(\rmd x, t-s, x-y)\rmd y,
    \\
\sP_1\left\{\zeta_t^M\in\rmd s,\;M_t\in\rmd y,\;\Gamma(t)\in\rmd x\right\}&
    =\delta_0(\rmd s)\delta_0(\rmd y)G^-(\rmd x, t, x)\1_{\{x<0\}}  \\
    +&\frac{\la_1}{\xi_0(t, x)}\left(x\sI_0(t, x)-\gamma_1x\xi_1(t, x)\sI_1(t, x)\right)
    \theta(t, x)\1_{\{0<x<\gamma_0t\}}\frac{\rmd x}{\gamma_0}\\
    &+\la_0F_1(\rmd s, s, y)G^-(\rmd x, t-s, x-y)\rmd y.
\end{aligned}
\]

\bibliographystyle{elsarticle-harv}

\begin{thebibliography}{99}
%
\bibitem{Luisa}
Beghin L., Nieddu L., Orsingher E., 2001.
Probabilistic analysis
of the telegrapher's process with drift by mean of relativistic
transformations. J. Appl. Math. Stoch. Anal.   14, 11--25.
%
\bibitem{BR}
  Bogachev L.,  Ratanov N., 2011. Occupation time distributions for the
telegraph process, Stochastic Processes and their Applications, 
121(8), 1816--1844. 

doi:10.1016/j.spa.2011.03.016 
%
\bibitem{cinque}
Cinque F. and Orsingher E., 2020.
On the  distribution of the maximum of the telegraph process 

\hskip-4mm Cinque F. and Orsingher E., 2020.
On the exact distributions of the maximum of the
asymmetric telegraph process.

\hskip-4mm Cinque F., 2020.
The negative reflection principle and the joint
distribution of the telegraph process and its
maximum. 

arXiv:2003.04044v1 9 Mar 2020

arXiv:2010.02689v1 6 Oct 2020

arXiv:2011.00342v1 31 Oct 2020
%
\bibitem{DiC2001}
Di Crescenzo, A. 2001.
On random motions with velocities alternating at Erlang-
distributed random times. Adv. Appl. Prob. 33, 690--701.
%
\bibitem{foong}
Foong S.K., Kanno S., 1994.
Properties of the telegrapherÕs random process with or without a trap, 
Stochastic Processes and their Applications, 53, 147--173.

doi: 10.1016/0304-4149(94)90061-2
%
\bibitem{GR}
Gradshteyn I.~S.,   Ryzhik, I.~M. 1994.
Table of Integrals,
Series, and Products, 5th ed.
Academic Press,
Boston.
%
\bibitem{KR}
Kolesnik A.D., Ratanov N., 2013.
Telegraph Processes and Option Pricing, 
Springer-Verlag, Heidelberg-New York-Dordrecht-London.

doi: 10.1007/978-3-642-40526-6
%
\bibitem{KS}
Karatzas, I.,   Shreve, S. 1998.     Methods of Mathematical Finance, 
Springer, Berlin.
%
\bibitem{JAP51}
 L\'opez O. and Ratanov N., 2014.
On the asymmetric telegraph processes.  {Journal of Applied Probability}, 
{51}(2), 569--589.

doi: 10.1239/jap/1402578644
%
\bibitem{O90}
Orsingher E. 1990.
\newblock
Probability law, flow function,
 maximum distribution of wave-governed random motions
and their connections with Kirchhoff's laws. 
\newblock
{Stoch. Process. Appl.} 34:49--66.
%
\bibitem{pinsky}
Pinsky M.A., 1991. 
Lectures on Random Evolution, World Scientific.

doi: 10.1142/1328 
%
\bibitem{RevuzYor}
Revuz D. and Yor M., 1999. 
Continuous Martingales and Brownian Motion (2nd ed.). Springer-Verlag, 
Berlin-Heidelberg-New York

doi: 10.1007/978-3-662-06400-9
%
\bibitem{SZ2004}
Stadje W. and Zacks S., 2004.
Telegraph processes with random velocities.
J. Appl.  Probab. 41(3),  665--678
%
\bibitem{zacks}
Zacks S., 2017. 
Sample Path Analysis and Distributions of Boundary Crossing Times,
Lecture Notes in Mathematics 2203, 
Springer International Publishing AG, 2017.

doi: 10.1007/978-3-319-67059-1
%
\end{thebibliography}

\end{document}